\documentclass{amsart}

\usepackage{amssymb,amsmath,amsthm}
\usepackage{color}
\usepackage{enumerate}
\usepackage{fullpage}
\usepackage{url}
\usepackage[unicode,pdfusetitle]{hyperref}
\usepackage[mathscr]{euscript}
\usepackage{graphicx}
\usepackage{mathtools}

\newtheorem{maintheorem}{Theorem}

\newtheorem{theorem}{Theorem}
\newtheorem*{theorem*}{Theorem}
\newtheorem{fact}[theorem]{Fact}
\newtheorem*{fact*}{Fact}

\newtheorem*{claim*}{Claim}

\newtheorem*{proposition*}{Proposition}
\newtheorem{lemma}[theorem]{Lemma}
\newtheorem*{lemma*}{Lemma}

\newtheorem*{question*}{Question}

\newtheorem*{conjecture*}{Conjecture}
\newtheorem{corollary}[theorem]{Corollary}
\newtheorem*{corollary*}{Corollary}

\theoremstyle{definition}
\newtheorem{definition}[theorem]{Definition}
\newtheorem*{definition*}{Definition}
\theoremstyle{remark}
\newtheorem{remark}[theorem]{Remark}
\newtheorem*{remark*}{Remark}

\newtheorem*{example*}{Example}

\numberwithin{theorem}{section}
\numberwithin{claim}{section}
\numberwithin{equation}{section}

\newcommand{\bbZ}{\mathbb{Z}}

\newcommand{\cL} {\mathcal{L}}
\newcommand{\bbR}{\mathbb{R}}

\newcommand{\bbC}{\mathbb{C}}\newcommand{\cR}{\mathcal{R}} \newcommand{\gen}[1]{\langle #1 \rangle}

\newcommand{\fj}{\mathfrak{j}}

\newcommand{\fg}{\mathfrak{g}}

\newcommand{\rexp}{\bbR_{\exp}}
\DeclareMathOperator{\trdeg}{trdeg}

\title{The Klein j-Function is not Pfaffian over the Real Exponential Field}
\author{Christoph Kesting}
\date{\today}
\email{kestingc@mcmaster.ca}
\address{Department of Mathematics and Statistics, McMaster University, Hamilton, Ontario, Canada}
\date{\today}

\begin{document}

\begin{abstract}
James Freitag showed that the  Klein $j$-function is not paffian over the complex numbers. We expand on this result by showing that a restriction of the Klein $j$-function to the imaginary interval $(0,i)$ is not pfaffian over the real field exponential field in the sense of  Miller and Speissegger. 
\end{abstract}
\maketitle
\section{Introduction}

Pfaffian functions are a class of functions introduced by Khovanskii \cite[Section 2.3]{khovanskiiFewnomials1991}, given by differential conditions, which have strong finiteness properties. These properties have played a crucial in showing the model completeness of certain o-minimal expansions of the reals. Verifying whether or not a function is pfaffian, or satisfies any of the generalizations of that concept remains a difficult problem. 
It is well known that any algebraic function is Pfaffian on a suitable domain and that differentially transcendental functions in one variable or those that violate the finiteness properties of the class are not Pfaffian. 
In \cite{freitagNotPfaffian2021} James Freitag introduced a new method to show that a function is not pfaffian. The key observation is that pfaffian functions satisfy by definition certain order one differential equations, which solutions to higher order strongly minimal differential equations can not satisfy. This results in a strong criterion for non-pfaffianness over $\bbC$. We aim to extend this result to the o-minimal setting with both the real field $\bar \bbR= (\bbR; + ,-, \times,0,1,<)$ and the real exponential field $\rexp=(\bbR; + , - , \times, 0,1, <, \exp)$. 

For this we expand the definition of pfaffian to a general o-minimal expansion $\cR$ of $\bar \bbR$ following \cite{miller2002pfaffian}.
\begin{definition}
An $\cR$-pfaffian chain is list of analytic functions $f_1,\dots,f_n:\bbR^m \to \bbR$ with $\cR$-definable functions $p_{i,j}:\bbR^{(i+1)m} \to \bbR^m $ such that for all $x \in \bbR^m$ we have:
$$ \frac{\partial f_i}{\partial x_j}(x) = p_{i,j}(x,f_1(x),\dots,f_i(x)).$$ 
We call a pfaffian chain classical if all the $p_{i,j}$ are polynomial as in \cite[Section 2.3]{khovanskiiFewnomials1991}.
\end{definition}
For the purposes of this paper, we will be restricting to pfaffian chains in one variable.

As in Freitags work, we choose to explore the results using the example of the Klein $j$-function. The $j$-function is an  important modular function in number theory, but we are using it as a stand-in solely for the differential-algebraic properties $j$ has.

The $j$-function satisfies an order $3$ algebraic differential equation, namely 

$$\left( \frac{y''}{y'} \right)'-\frac{1}{2}\left( \frac{y''}{y'}\right)^2 + (y')^2 \cdot \frac{y^2 -1968y+2654208}{y^2 (y-1728)}=0.$$

As conjectured by Mahler, Nishioka could show that the $j$ function does not satisfy any lower-order differential equation:
\begin{fact}[\cite{ConjectureMahlerAutomorphic1989}]
    The $j$-function satisfies no non-trivial differential equation of degree $\leq 2$.
\end{fact}

Using this fact and stability theoretic machinery, Tom Scanlon and James Freitag established the following result:
\begin{fact}[\cite{freitagStrongMinimalityJfunction2014}]
    The differential equation the $j$-function satisfies is strongly minimal in $DCF_0$.
\end{fact}

\begin{remark}\label{rem:strong minimality}
In general, if a function $f$ satisfies a strongly minimal differential equation of order $n$ over $K$, then the differential equation is irreducible over $K$ and, if $F \supseteq K$ is a differential field extension, then the transcendence degree of the field extension $F(f,f', \dots ,f^{(n)})$ is 
$$\trdeg(F(f,f', \dots ,f^{(n)})/F) = 0 \text{ or } n,$$
depending on whether or not $f\in F$.
\end{remark}
For the  example of the $j$-function this means:

\begin{corollary} \label{cor:j extensions}
    If $F \supseteq \bbC$ is a differential field extension, then the transcendence degree of the field extension $F(j,j',j'')$ is
$$\trdeg(F(j,j',j'')/F) = 0 \text{ or } 3,$$
depending on whether or not $j\in F$. 
\end{corollary}

A recent applications of pfaffian functions emerged within the context of Pila-Wilkie counting \cite{pila2006rational}. In particular, in \cite{binyamini2023effective}, an effective version of Pila-Wilkie counting for sets existentially definable from pfaffian functions was shown. This found applications in the same paper in an effective and uniform version of the Manin-Mumford conjecture, which is a statment about limiting torsion points of abelian varties. 
The original point counting result has been instrumental in the proof of the André–Oort conjecture in \cite{pila2011minimality} and see  \cite{scanlon2017minimality} for a survey of the area. These proofs  are using the o-minimal structure $\bbR_{an,\exp}$ with function symbols for all real analytic function restricted  to $[0,1]$ and the full real $\exp$, mainly for the definability for a the restriction of the  $j$-function. The reliance on $\bbR_{an,\exp}$ in the proofs precludes an effective version of Pila-Wilkie counting to happen, and thus an  effective version of the André–Oort conjecture would require a different structure or counting argument, as for example in  \cite{binyamini2024log}.
Beyond that, in the respective subsections proving the theorems, we also go into detail to what extent these results generalize beyond the $j$-function.

Below, we will be working with $\fj=j(i \cdot x): (0,1) \to  \bbR$.
Furthermore, for an o-minimal field $\cR$ and some element $a$, we denote by $\cR \gen{a}$ the o-minimal structure generated by taking the definable closure of $\cR$ and $a$.

Using this characterization, we will prove the following two non-pfaffianess results.

\begin{maintheorem}\label{thm:mainA}
 $\fj$ is not algebraic over any  $\bar \bbR$-pfaffian chain.
\end{maintheorem}
By this we mean  that a germ of $\fj$ is not algebraic over the germs of any $\bar \bbR$-pfaffian chain. In the proof this reduces to generated o-minimal $\cL_\text{ring}$-substructures of the $\cL_\text{ring}$  reduct of the  Hardy field of the pfaffian closure $\bbR_\text{pfaff}$, as $\fj$ lives in  $\bbR_\text{pfaff}$\footnote{upcoming paper of Gareth Jones, Patrick Speissegger and Jonathan Kirby} and algebraicity is $\cL_{ring}$-definability in this setting.
Using an analogous model-theoretic description  for $\bbR_{\exp}$, we can consider the $\cL_{\exp}$-closure of the germs of any pfaffian chain in the $\cL_{\exp}$-reduct of said Hardy field and prove: 
\begin{maintheorem}\label{thm:mainB}
     $\fj$ is not $\cL_{\exp}$-definable over any $\rexp$-pfaffian chain.
\end{maintheorem}

In the respective sections proving the theorems, we also go into detail to what extent these results generalize beyond the $j$-function.
But first, we will need the following lemma to control field extensions by $\cR$-paffian chains.

\begin{lemma}
\label{minimal pfaffian chains}
 Let $f_1,\dots,f_n$ be an $\cR$-pfaffian chain. Then there exists a minimal length $\cR$-paffian chain $g_1,\dots,g_m$ such that $\cR \gen{f_1,\dots,f_n}=\cR \gen{g_1,\dots,g_m}$. In particular for $\cR= \bar\bbR$ we have that the real algebraic closure of $\bbR (g_1 ,\dots, g_m)$ is a differential field and that  $\trdeg(\bbR (g_1 ,\dots, g_m)/\bbR)=m$.
\begin{proof}
Suppose not.
Let $f_1,\dots,f_n$ with corresponding functions $p_{i,j}$ for the derivatives  be given. Let $k\leq n$ be maximal such that $(f_j)$ for $j<k$ are algebraically independent over $\bbR$. Then $f_k$ is algebraic over $f_1,\dots,f_{k-1}$. 
By definable choice, $f_k$ is already definable from $f_{1},\dots,f_{k-1}$, using some definable function $Q(x,f_1,\dots,f_{k-1})$.
Now we can remove $f_k$ from the chain by substituting all instances of $f_k$ in $p_{i,j}$ for $i\geq k$ with $Q$.
Then $f_1,\dots,f_{k-1},f_{k+1},f_n$ is a $\cR$-pfaffian chain of lenght $n-1$, generating the same extension of $\cR$ as the chain we started with.
\end{proof}
\end{lemma}
\subsection*{Acknowledgements}
Special thanks go out to Gareth Jones for his helpful suggestions for approaching this problem and James Freitag for further advice and feedback. Also, I'd like to thank my PhD supervisor Patrick Speissegger for the countless hours spend discussing this problem. 

\section{The semialgebraic case}
In this section, we are working in the o-minimal structure $\bar \bbR= (\bbR, + ,- ,\times,0,1,<)$. Every set definable in this structure is semi-algebraic. As $\bar \bbR$-pfaffian functions are analytic by definition, we will be using the following fact, to generate algebraic dependences.

\begin{fact}[{\cite[Proposition 8.1.8]{bochnakRealAlgebraicGeometry1998}}]\label{fact:analyticsemialgebraic}
On an open semialgebraic subset $U\subset \bbR^n$, for a function $f:U \to \bbR$ the following is equivalent:
\begin{itemize}
    \item $f$ is analytic and semi-algebraic on $U$;
    \item $f$ is analytic and algebraic on $U$.
\end{itemize}
\end{fact}

With  this, we can  go on to prove our first result:

\begin{theorem}[Theorem A]  $\fj$ is not algebraic over any  $\bar \bbR$-pfaffian chain..
\begin{proof}
Suppose it is. Let $f_{1},\dots,f_{n+1}$ be a  $\bar \bbR$-pfaffian chain such that  $\fj$ is definable over $f_1,\dots,f_{n+1}$ and $f_1,\dots,f_{n+1}$ is of minimal length and minimal in the sense of \ref{minimal pfaffian chains}.
Then by minimality, we may assume that $\fj$ is transcendental over $f_1,\dots,f_n$ and inter-definable with $f_{n+1}$.
By o-minmality, this is witnessed by definable functions $G$ and $H$ such that $j=G(f_1,\dots,f_{n+1})$ and $f_{n+1}=H(f_1,\dots,f_n,j)$. 
Then let $X\subseteq \bbR^{n+3}$ be a semialgebraic set with $(f_{1}(t),\dots,f_{n}(t),\fj(t),\fj'(t),\fj''(t)) \in X$ for $t\in (0,1)$.
As, by the pfaffian condition and chain rule, $\fj'(t)$ and $\fj''(t)$ are definable from the first $n+1$ variables of the tuple over $\bbR$, we may intersect $X$ with the graph of the functions defining them. 
Then $X$ is of dimension $n+1$ and we may assume that $X$ is given as the graph of semialgebraic function $F$. Without loss of generality, we assume that $F$ satisfies  $$F(f_{1}(t),\dots,f_{n}(t),\fj(t))=(\fj'(t)),\fj''(t))\text{ for } t \in (0,1).$$
Now by analytic cell decomposition, $F$ and  is analytic on an open set $U$.
Then again by Fact \ref{fact:analyticsemialgebraic}, we have that $F$ satisfies $$p(F( x_{1},\dots,x_{n+1}),x_{1},\dots, x_{n+1})=0$$ for a polynomial $p$ on $U$. That means that $F$ is an algebraic function on $U$, so $X$ is an algebraic set. Over $(\bbR,+,-,\times,0,1)$ transcendence degree and semialgebraic dimension coincide. In particular $\dim(X)=\trdeg(X)=n+1$, where $\trdeg(X)$ is the maximal transcendence degree of a point in $X$ over $\bbR$.  Now take a generic point  $p \in (0,1)$ such that $(f_1(p),\dots,f_n,\fj(p))$ is generic in $U$. 
Next, we define $\bbR\gen{f_{1},\dots,f_{n},\fj,\fj',\fj''}$  as the function field generated by the germs of $f_{1},\dots,f_{n},\fj,\fj',\fj''$ at $p$.  
Then, by the minimality of the $\bar \bbR$-pfaffian chain for all $i\leq n$ we have: 
$$ \trdeg(\bbR\gen{f_{1},\dots,f_{i}}/\bbR \gen{f_{1},\dots,f_{i-1}})=1$$
so $\trdeg(\bbR \gen{f_{1}\dots,f_{n}}/\bbR)=n$ and $\trdeg(\bbR \gen{f_{1}\dots,f_{n},\fj,\fj',\fj''}/\bbR \gen{f_{1}\dots,f_{n}})={1}$. By the identity theorem for analytic functions we have that $\trdeg(\bbC \gen{\tilde f_{1}\dots, \tilde  f_{n},\tilde \fj,\tilde \fj',\tilde \fj''}/\bbC \gen{\tilde f_{1}\dots,\tilde f_{n}})={1}$, where $\tilde h$ for a real analytic germ $h$ is its complex continuation. This contradicts the strong minimality of the differential equation of $\fj$ as characterized in \ref{cor:j extensions}.
\end{proof}
\end{theorem}

\begin{corollary} \label{cor:notrealpfaffian}
 $\fj$ is not $\bar\bbR$-pfaffian.
\begin{proof}
 $\fj$ would be algebraic over any such $\bar \bbR$-pfaffian chain, contradicting the theorem.
\end{proof}
\end{corollary}

\begin{corollary}
    No real reduct of a  function satisfying a strongly minimal algebraic differential equation over $\bbC$ of order $\geq 2$ is algebraic over an $\bar \bbR$-pfaffian chain.
    \begin{proof}
        This captures precisely the properties of the $j$-function that go into the proof of the theorem.
    \end{proof}
\end{corollary}

\begin{corollary}\label{cor:derivative}
    Let $g$ satisfy a strongly minimal algebraic differential equation over $\bbC$ of order $m$. Then for all $k\leq m-2$, any real restriction of the $k$-th derivative of $g$ is not algebraic over any $\bar \bbR$-pfaffian chain.
    \begin{proof}
        Assume there is a minimal  $\bar \bbR$-pfaffian chain $f_1,\dots,f_n$  of length $n$ such that a restriction $\fg^{(k)}:(0,1) \to \bbR$ of $g^{(k)}$ with $k\leq m-2$ is algebraic over $f_1,\dots,f_n$. Then we can extend the pfaffian chain to $f_1,\dots,f_n,\fg^{(k)},\fg^{(k-1)},\dots,\fg',\fg$. Then there  exists a semialgebraic set  $X$  with $(f_1,\dots,f_n,\fg^{(i)},\fg^{(i-1)},\dots,\fg',\fg)\in X$ of dimension $n+k$. Following the proof of the Theorem verbatim, as  $n+k<n+m-1$, we again get a transcendence degree inequality contradicting strong minimality. 
    \end{proof} 
\end{corollary}

\section{The real exponential case}

To apply a similar strategy as in the semialgebraic case to the exponential case, we have to establish some facts.
First, we need to have an appropriate analytic cell decomposition result in $\rexp$, and we need the right Ax-Schanuel type statement that connects the transcendence degree of $\exp$ and the $j$-function.

For the first part, we need the following results coming from work in o-minimality:

\begin{fact}[\cite{wilkieMODELCOMPLETENESSRESULTS}]\label{thm:modelcompleteness}
$\bbR_{\exp}$ is model complete i.e. every first-order definable set is equivalent to a projection of a quantifier-free definable set. 
\end{fact}

\begin{fact}[Analytic cell decomposition for $\bbR_{\exp}$ \cite{vandendriesRealExponentialField1994}] \label{thm:analytic cell decomposition}
The real exponential field $\rexp$ has analytic cell decomposition i.e. for every finite collection of definable subsets $A_i$ of $\bbR^n$ there exists a partition of $\bbR$ into special well-behaved definable sets $C_j$ called analytic cells. These cells respect the $A_i$ in the sense that $A_i \cap C_j = C_j$ or $\emptyset$. For more details on the analytic side, see \cite[Section 8]{vandendriesRealExponentialField1994}.
\end{fact}

\begin{lemma}
Let $f_{1},\dots,f_{n}$ be an $\bbR_{\exp}$-pfaffian chain. Then for each $i\leq n$ there are $\bar \bbR$-pfaffian chains $g_{i,1},\dots ,g_{i,m_{i}}$ in $n_{i}$ variables such that the real algebraic closure of $$\bbR(f_{1},g_{1,1}(f_{1}),\dots ,g_{1,m_{i}}(f_{1}),f_{2},g_{2,1}(f_{1}, f_{2}), \dots , g_{i-1,m_{i-1}}(f_{1}, \dots, f_{i-1}),f_{i},g_{i,1}(f_{1}, \dots, f_{i}),\dots ,g_{i,m_{i}}(f_{1}, \dots, f_{i}))$$ is a differential field extension  and the transcendence degree of the extension over $\bbR$ matches the number of generators.
\begin{proof}By induction assume this already holds for the first $i$ elements of the chain and let $\cR$ be the resulting differential field.  Then let $\partial f_{i+1}/\partial x = h(x,f_{1}(x),\dots,f_{n+1}(x))$ with $h$ being $\bbR_{\exp}$-definable. Then by analytic cell decomposition and model completeness for $\Gamma_h$ the graph of $h$ we have on a cell $C$ that
$$\bbR_{\exp}\models (x,w,y)\in \Gamma_{h} \Leftrightarrow \exists \bar zF(x,w,y,z)=0$$ for a quantifier free $\bbR_{\exp}$-definable $F$. This $F$ is now however given as $p(x,w,y,z,g_{1}(x,w,y,z), \dots,g_{m}(x,w,y,z))$ where $p$ is a polynomial and $g_{1},\dots,g_{m}$ is a $\bbR$-pfaffian chain of $\cL_{\exp}$-terms in $n_{i}$-variables.
Then subbing $f_1,\dots,f_{i+1}$ into $w$ yields the field extension: $$\cR(f_{i+1}(x),g_1(x,f_{1},\dots,f_{i+1},y,z),\dots g_m(x,f_{1},\dots,f_{i+1},y,z))$$ over which $\partial f_{i+1}/\partial x$ is $\bbR$-definable, hence semi-algebraic and therefore algebraic.
For transcendence degree matching the number of generators it is sufficient to follow the proof strategy of Lemma \ref{minimal pfaffian chains} to eliminate any redundancies.  
\end{proof}
\end{lemma}

As the second part of our story's  prerequisites, we need the following definition in our Ax-Schanuel statement. 
\begin{definition} \label{def:formal parametrization}
    A formal parametrization of a neighbourhood of a point $p\in \bbC$ in the sense of the following theorem is a power series of the form $t \in p + (\mathfrak{m} \setminus \{0\})$, where $\mathfrak{m}$ is the maximal ideal in $\bbC [[z_1,\dots,z_m]]$.
\end{definition}

\begin{fact}[Case $l=1$ of {\cite[Theorem E]{blázquezsanz2022differential}}]    \label{Ax-Schanuel}
Let $\hat p_{1},\dots,\hat p_{k+1}$ be formal parametrisations (in variables $z_{1},\dots,z_{m}$) of neighbourhoods of points $p_{1},\dots,p_{k+1}$ in $\mathbb{H}$. If the transcendence degree:
$$ \trdeg (\bbC(\hat p_{1},\dots, \hat p_{k+1},\exp (\hat p_{1}), \dots,\exp(\hat p_{k}),j(\hat p_{k+1}),\dots,j''(\hat p_{k+1}))/ \bbC )< k +3 +rank\begin{pmatrix} \frac{\partial \hat p_{j}}{\partial \hat z_{i}} \end{pmatrix} $$
then there exists $n \in \bbZ^{k}\setminus \{0\}$ such that $\sum_{i=1}^{k}n_{i}\hat p_{i}\in \bbC$.
\end{fact}

\begin{remark}
A more general version of the Ax-Schanuel result holds for fuchsian uniformizing functions of hyperbolic curves instead of $j$ in multiple variables, see \cite[Section 8]{blázquezsanz2022differential}. A similar result to the following theorem can be established as well.
\end{remark}

With our prerequisites out of the way, we can now prove Theorem B.

\begin{theorem}[Theorem B]
$\fj$ is not $\cL_{\exp}$-definable using any $\rexp$-pfaffian chain.
\begin{proof}
Assume not. Let $f_{1},\dots, f_{n+1}$ be an $\bbR_{\exp}$-pfaffian chain such that $\fj$ is definable over the chain. By the  lemma above, we can find   $g_{1,1},\dots g_{n+1,m_{n+1}}$ such that the resulting differential field has minimal transcendence degree and $f_{n+1}$ is inter-definable with $\fj$ over the rest of the generators, of which there are $\sum_{i=1}^{n}n_{i}+ n=:M$ many.
By o-minimality, this is witnessed by definable functions $G$ and $H$ such that $\fj=G(t,f_1,g_{1,1},\dots,f_{n+1})$ and $f_{n+1}=H(t,f_1,g_{1,1},\dots,f_n,\fj)$. 
Then there exists a definable set $X \subseteq \bbR^{M+4}$  such that for $t \in (0,1)$ we have $$(t,f_{1}(t),g_{1,1},\dots,g_{n,n_{i}}(t),\fj(t),\fj'(t),\fj''(t))\in X.$$
By the pfaffianness and the chain rule, we can intersect $X$ with the graphs of the functions defining $\fj'=\frac{\partial}{\partial t} G(f_1 ,g_{1,1},\dots,f_n,H(f_1,g_{1,1},\dots,f_n,\fj))$ and $\fj''= \frac{\partial}{\partial t} \fj'$ analogously.
Then by analytic cell decomposition \ref{thm:analytic cell decomposition}, $X$ is given on a Cell $C$ of maximal dimension in $X$ as the graph of a function $h$, without loss of generality satisfying  $$h(t,f_{1}(t),g_{1,1},\dots,f_{n},\dots, g_{n,n_{i}}(t),\fj(t))=(\fj'(t)),\fj''(t)) \text{ for } t \in (a,b),$$
for an open interval $(a,b)\subseteq (0,1)$ depending on the pfaffian chain and the cell $C$.
Then by model completeness \ref{thm:modelcompleteness}, we have that 
$$ \bbR_{\exp}\models (t,w,x,y) \in \Gamma_{h}\Leftrightarrow \exists z F(t,w,x,y,z)=0 ; \quad \text{ with }|z|=N $$
for $F : \bbR^{1+M+3+N} \to  \bbR$ quantifier free definable as a composition of a polynomial with $g_{1},\dots,g_{l}$  a $\bar \bbR$-pfaffian chain of $\cL_{\exp}$-terms. rearranging all $\cL_{\exp}$-terms from the new $g_{1},\dots,g_{l}$ potentially increased $N$ yields a a polynomial $Q$ such that $$F=Q(t,w,x,y,z,\exp t, \exp w, \exp x, \exp y, \exp z).$$Now we take a generic point  $p \in (0,1)$ such that $(f_{1}(p),\dots,g_{n,n_{i}},\fj(p),\fj'(p),\fj''(p))$ is generic in $X\cap C$.  

Subbing $f_{1},g_{1,1},\dots,f_{n},\dots, g_{n,n_{i}}$ into $w$, $\fj(t)$ into $x$ and $(\fj',\fj'')$ into $y$ we get the following field extension: 
$\bbR(t,f_{1},g_{1,1},\dots, g_{n,n_{i}},\fj(t),\fj'(t),\fj''(t),\exp(t),\exp(f_{1}),\exp(g_{1,1}),\dots,\exp(g_{n,n_{i}}),\exp(\fj(t)),\exp(\fj'(t)),\exp(\fj''(t)))$ 
By adding additional variables $p_{0},\dots,p_{3}$ for the power series $t,\fj(t),\fj'(t),\fj''(t)$ respectively we can bring this in a form that is compatible with the Ax-Schanuel Theorem \ref{Ax-Schanuel}. Note that these $p_0,\dots,p_3$ can not contribute to the transcendence degree of the extension, $\fj(t),\dots,\fj''(t)$ can increase the transcendence degree by at most $1$ by pfaffianess, there are $M$ many $f_i$ and $g_{kl}$, and $z$ is an $N$ tuple. So, we have the following bound for the transcendence degree:
$$\begin{aligned}
&\trdeg(\bbR( t,p_{0},\dots,p_{3},f_{1},g_{1,1},\dots,, g_{n,n_{i}},\fj(t),\fj'(t),\fj''(t), z, \\ &\exp(p_{0}), \dots, \exp(p_{3}),\exp(f_{1}),\exp(g_{1,1}),\dots,\exp(g_{n,n_{i}}), \exp(z)/\bbR) \\ 
\leq & 1 + M +1+ N +4+ M  +N \\
= &2 (M+N)+6\\
< &2 (M+N)  +7\\
\leq & \underbrace{M  +N + 4}_{\text{variables in}\exp} + \underbrace3_{\text{variables in }\fj,\fj',\fj''} + \underbrace{rank\left( \frac{\partial t,\partial p_{i},\partial f_{i}, \partial g_{i,l}, \partial z}{\partial s_{j}}\right)}_{\geq M  +N, \text{ by minimality of the chain and $z$ being independent}}
\end{aligned}
$$
where $s_j$  are the variables that we are working in.
Thus   taking the complex analytic continuations and applying Ax-Schanuel \ref{Ax-Schanuel} we have that there is an non zero linear integer combination of the variables plugged into $\exp$ that is constant. In particular we have that $$p_{0},\dots,p_{3},f_{1},g_{1,1},\dots, g_{n,n_{i}}, z$$are algebraic over $\bbC$. But subbing back in $t,\fj(t),\fj'(t),\fj''(t)$ for $p_{0},\dots,p_{3}$ we can then use the minimality of chain expansion to show that the differential field extension 
$$\trdeg(\bbR (f_{1},g_{1,1},\dots, g_{n,n_{i}})/\bbR)=M.$$
Then the expansion by both $z$ and $t$ remains a differetial field 
with $$\trdeg(\bbR(f_{1},g_{1,1},\dots, g_{n,n_{i}},t, z)/\bbR)=M+N+1$$
still witnessing no dependencies. Now by strong minimality of $j$ we have that either $$\trdeg(\bbR(f_{1},g_{1,1},\dots, g_{n,n_{i}},t, z, \fj(t),\fj'(t),\fj''(t))/\bbR)$$
is $M+N+1$ or $M+N+4$ with the former violating the minimality of the pfaffian chain and the latter violating the dependence from Ax-Schanuel. 
In either case we have a contradiction. 
\end{proof}
\end{theorem}

\begin{corollary}
   $\fj$ is not $\bbR_{\exp}$-pfaffian.
    \begin{proof}
        If $\fj$ was $\rexp$-pffaffian, it would be definable over it's own pfaffian chain.
    \end{proof}
\end{corollary}

\begin{remark}
    The other corollaries from the semialgebraic case do not translate as easily to other solutions of strongly minimal differential equations as the corresponding Ax-Schanuel result is oftentimes not known.  
\end{remark}

\end{document}